\documentclass[12pt]{article}
\usepackage{amsmath, amssymb, eucal, latexsym}
\usepackage[cp1251]{inputenc}
\usepackage[english]{babel}
\usepackage{graphicx}

\def\Z{\Bbb Z}
\def\N{\Bbb N}
\def\R{\Bbb R}
\def\ve{\varepsilon}
\def\oT{{\mathbf T}}
\def\v{\vec}
\def\D{\Delta}
\def\ox{{\mathbf x}}
\def\ow{{\mathbf w}}
\def\mes{\operatorname{mes}}

\newtheorem{theorem}{Theorem}

\newtheorem{lemma}[theorem]{Lemma}

\newtheorem{note}[theorem]{Remark}

\begin{document}

 \begin{center}
 \textbf{ON WIENER NORM OF SUBSETS OF $\Z_p$ OF MEDIUM SIZE}
 \end{center}

 \begin{center}
                                                         S. V. KONYAGIN
\footnote{
The first author is supported by grant RFBR 14-01-00332
and grant Leading Scientific Schools N 3082.2014.1}                                                         ,
I. D. SHKREDOV\footnote{
The second author is supported by grant
mol\underline{ }a\underline{ }ved 12--01--33080.
}\\

    \end{center}

\bigskip

\begin{center}
    Abstract.
\end{center}

{\it \small
    We
    give a lower bound
    for  Wiener norm of characteristic function of
subsets $A$ from $\Bbb Z_p$, $p$ is a prime number, in the situation when
$\exp\left((\log p/\log\log p)^{1/3}\right) \le |A| \le p/3$.}

\bigskip
\section{Introduction}
\bigskip

We consider the abelian group $G=\Z_p=\Z/p\Z$, where $p$ is a prime number.
Denote the Fourier transform of a complex function on $G$ to be a new function
$$\hat f(\gamma)=\frac1p\sum_{x\in G}f(x)e_p(x\gamma) \,,$$
where $e_p(u)=\exp(2\pi iu/p)$ (we note that $e_p$ is correctly defined for
$u\in\Z_p$). It is known that the function $f$ can be reconstructed from
$\hat f$ by the inverse Fourier transform
\begin{equation}\label{Fourier_inv}
f(x)=\sum_{\gamma\in\Z_p}\hat f(\gamma)e_p(-x\gamma).
\end{equation}
We define the Wiener norm of a function $f$ as
$$\|f\|_{A(G)}=\|f\|_A=\|\hat f\|_1=\sum_{\gamma\in\Z_p}|\hat f(\gamma)| \,.$$
By $\chi_S$, $S\subset G$ denote the characteristic function of some set $S$.

In this note we discuss the problem of estimation from below the Wiener norm
of $\chi_A$ for $A\subset\Z_p$ in terms of $p$ and $|A|$.

If $x\in A$, then, by~(\ref{Fourier_inv}), we have
$$1=\left|\sum_{\gamma\in\Z_p}\hat f(\gamma)e_p(-x\gamma)\right|
\ge\sum_{\gamma\in\Z_p}|\hat f(\gamma)| \,.$$ Thus, we get a trivial
estimate for Wiener norm of any nonempty $A\subset\Z_p$
\begin{equation}\label{trivial}
\|\chi_A\|_{A}\ge1.
\end{equation}
Next we observe that because of
$$\|\chi_{\Z_p\setminus A}\|_{A}=\|\chi_A\|_{A}+(1-2|A|/p)$$
it is sufficient to consider the case $|A|<p/2$.
It is easy to see that if $A\subset\Z_p$ is an arithmetic progression with
\begin{equation}\label{nat_A}
2\le|A|<p/2
\end{equation}
then
$$\|\chi_A\|_A\asymp\log|A|.$$
It is commonly believed that for any $A$ satisfying~(\ref{nat_A})
there is the same lower bound
\begin{equation}\label{conj_lower_bound}
\|\chi_A\|_A\gg\log|A|.
\end{equation}

The first nontrivial lower bound for $\|\chi_A\|_A, |A|<p/2$, in some range was
established in\cite{GK}:
$$\|\chi_A\|_A\gg\frac{|A|}{p}\left(\frac{\log p}{\log\log p}\right)^{1/3}.$$
This estimate was improved by T.~Sanders \cite{Sanders} for $|A|<p/2$, $|A|\gg p$.
As was shown in \cite{KS}, the results of \cite{Sanders} imply the following.

\begin{theorem}
\label{Char}
Let $p$ be a prime number, $A\subset\Z_p$, $0<\eta=|A|/p<1/2$. If $\eta\ge(\log p)^{-1/4}(\log\log p)^{1/2}$ then
$$\|\chi_A\|_{A}\gg(\log p)^{1/2}(\log\log p)^{-1}
\eta^{3/2} \left(1 + \log\left(\eta^2(\log p)^{1/2}(\log\log p)^{-1}\right)\right)^{-1/2},$$
and if  $\eta<(\log p)^{-1/4}(\log\log p)^{1/2}$ then
$$\|\chi_A\|_{A}\gg\eta^{1/2}(\log p)^{1/4}(\log\log p)^{-1/2}.$$
\end{theorem}
Our interest to study Wiener norm of large subsets of $\Z_p$
was inspired by the paper of V.V.~Lebedev \cite{Lebedev} on
quantitative variants of Beurling--Helson theorem.

Theorem \ref{Char} is nontrivial if our subset $A$ is large, that is
$$|A|p^{-1}(\log p)^{1/2} (\log\log p)^{-1} \to\infty $$
(and of course $|A|<p/2$). For small $A$ we proved in \cite{KS} a sharp estimate.
\begin{theorem}
\label{Charsmall}
Let $p$ be a prime number, $A\subset\Z_p$, and
$$2\le|A|\le\exp\left((\log p/\log\log p)^{1/3}\right).$$
Then
$$\|\chi_A\|_{A}\gg\log |A|.$$
\end{theorem}

In this note we study the subsets $A\subset\Z_p$ of medium size.
Our main result is the following assertion.
\begin{theorem}
\label{mediumsize}
Let $p$ be a prime number, $A\subset\Z_p$,
$$\exp\left((\log p/\log\log p)^{1/3}\right)\le|A|\le p/3.$$
Then
$$\|\chi_A\|_{A}\gg(\log(p/|A|))^{1/3}(\log\log(p/|A|))^{-1+o(1)}.$$
\end{theorem}

\bigskip

We observe that using arguments of Theorem~\ref{Charsmall} one can get
analogious estimates
for sets $A$ slightly exceeding the bound indicated in the statement.
However, the improvement is marginal. Moreover, it seems that by that way
one cannot get a nontrivial estimate for rather large subsets,
namely, such that $\log|A|\gg\log p$.

\bigskip
\section{Comparison with the continuous case}
\bigskip

We denote $e(u)=\exp(2\pi iu)$.
For sets $B\subset\Z$ a continuous analog of (\ref{conj_lower_bound})
is a well--known fact. Namely, it was proved in \cite{Kon}
and \cite{MGPS} that if $B\subset\Z$, $2\le|B|<\infty$ then
$$\int_0^1\left|\sum_{b\in B}e(bu)\right|du\gg\log|B|.$$
Moreover, in \cite{MGPS} the following stronger result was proved:
if $b_1<\dots<b_l$ are real numbers and $c_j$ are arbitrary complex numbers then
\begin{equation}\label{Hardy}
\int_0^1\left|\sum_{j=1}^l c_je(b_ju)\right|du\gg\sum_{j=1}^l\frac{|c_j|}j.
\end{equation}
This inequality implies the following lemma.
\begin{lemma} \label{n,2n}
Let $n\in\N$, $B\subset[-2n,2n] \subset\Z$, $|B|\ge2$, $0<\eta<1/2$,
$|B\cap[-n,n]|\ge(1-\eta)|B|$, $c(b)\,(b\in B)$ are complex numbers
with $c(b)=1$ for $b\in B\cap[-n,n]$. Then
$$\int_0^1\left|\sum_{b\in B} c(b)e(bu)\right|du\gg
\min\left(\log\frac1{\eta}, \log|B|\right).$$
\end{lemma}

{\bf Proof} Let $B=\{b_1<\dots<b_l\}$ where $l=|B|$,
and let $B\cap[-n,n]=\{b_{l_1}<\dots<b_{l_2}\}$.
The polynomial $\sum_{b\in B}^l c(b)e(bu)$ can be rewritten as
$\sum_{j=1}^l c_je(b_ju)$ where $c_j=1$ for $l_1\le j\le l_2$.
We denote
$$S=\int_0^1\left|\sum_{b\in B} c(b)e(bu)\right|du.$$
By (\ref{Hardy}),
$$S\gg\sum_{j=l_1}^{l_2}\frac 1j
\gg\log((l_2+1)/l_1).$$
We have $l_2-l_1+1\ge(1-\eta)l$. If $\eta<1/l$, then $l_1=1$, $l_2=l$,
$S\gg\log((l_2+1)/l_1)=\log l$ as required. If $\eta\ge1/l$, then we have
$$l_1\le\eta l+1<2\eta l.$$
Hence,
$$\log((l_2+1)/l_1)\ge\log((l_1+(1-\eta)l)/l_1)\ge\log((1+\eta)/2\eta)
\gg\log(1/\eta),$$
and we again get the assertion of the lemma.
$\hfill\Box$

\bigskip

The discrete and continuous $L^1$--norms of trigonometric polynomials
can be compared by the following lemma.
\begin{lemma} \label{disc_cont}
We have
$$\frac1p\sum_{\gamma\in\Z_p}\left|\sum_{|x|\le p/3}c_xe_p(x\gamma)\right|
\gg\int_0^1\left|\sum_{|x|\le p/3} c_xe(xu)\right|du.$$
\end{lemma}
See \cite{Zyg}, chapter 10, Theorem 7.28.

\bigskip

One can deduce (\ref{conj_lower_bound}) from Lemma~\ref{disc_cont}
provided that $A\subset[-p/3,p/3]$ (this inclusion means that any residue
$a\in A$ has an integer representative from $[-p/3,p/3]$) or
if some non--degenerate affine image of $A$ in $\Z_p$ is contained in
$[-p/3,p/3]$. This argument was used in the proof of Theorem~\ref{Charsmall}.

\bigskip

Now let us define the de la Vall\'ee-Poussin polynomials and means.
For functions
$$F(\gamma)=\sum_{x\in\Z_p}c_xe_p(x\gamma),\quad
G(\gamma)=\sum_{x\in\Z_p}d_xe_p(x\gamma)$$
we define their convolution
$$F*G(\gamma)=\sum_{x\in\Z_p}c_xd_xe_p(x\gamma).$$
It is easy to see that
$$F*G(\gamma)=\frac1p\sum_{\xi_1+\xi_2=\gamma}F(\xi_1)G(\xi_2).$$
Therefore,
\begin{equation}\label{L1_convol}
\sum_{\gamma\in\Z_p}|F*G(\gamma)|\le\frac1p\sum_{\gamma\in\Z_p}|F(\gamma)|
\sum_{\gamma\in\Z_p}|G(\gamma)| \,.
\end{equation}

Study of arbitrary trigonometric polynomials in $\Z_p$ can be
reduced to polynomials of small degree using de la
Vall\'ee-Poussin means. Define the de la Vall\'ee-Poussin
polynomial of order $n\le p/4$ as
$$V_n(\gamma)=\sum_{|x|\le n}e_p(x\gamma)
+ \sum_{n<|x|\le 2n}\frac{2n-|x|+1}{n+1}e_p(x\gamma)$$
and the de la Vall\'ee-Poussin mean for $F$ of order $n\le p/4$ as
$F*V_n$.

We need in the lemma.

\begin{lemma} \label{V_Pnorm}
For $n\le p/4$ the following inequality holds
$$\sum_{\gamma\in\Z_p}|V_n(\gamma)|\le3p.$$
\end{lemma}
The proof is contained in the proof of Theorem 7.28 of chapter 10
in \cite{Zyg}.

\bigskip

Using Lemma~\ref{V_Pnorm} and (\ref{L1_convol}) we
obtain
the following lemma.
\begin{lemma} \label{V_P}
For $n\le p/4$ the following inequality holds
$$\sum_{\gamma\in\Z_p}\left|\sum_{|x|\le n}c_xe_p(x\gamma)
+ \sum_{n<|x|\le 2n}\frac{2n-|x|+1}{n+1}c_xe_p(x\gamma)\right|
\le 3\sum_{\gamma\in\Z_p}\left|\sum_{|x|\le p/2}c_xe_p(x\gamma)\right|.$$
\end{lemma}

\bigskip

Combining Lemmas \ref{V_P}, \ref{disc_cont}, and \ref{n,2n} we get the following.
\begin{lemma} \label{trig}
Let $B\subset\Z_p$, $n\le p/6$, $0<\eta<1/2$. Assume that $|B\cap[-2n,2n]|\ge2$
and
$$
    |B\cap[-n,n]|\ge(1-\eta)|B\cap[-2n,2n]| \,.
$$
Then
$$\|\hat{\chi}_B \|_1\gg\min\left(\log\frac1{\eta}, \log|B\cap[-2n,2n]|\right).$$
\end{lemma}

\bigskip
\section{Balog-- Szemer\'edi-- Gowers theorem, Freiman's theorem, and
structure of sets with small Wiener norm}
\bigskip

Given an arbitrary set $Q\subset\Z_p$ and $k\in\N$, denote the quantity
$\oT_k (Q)$  as the number of solutions to the equation
$$x_1+\dots+x_k = x'_1+\dots+x'_k$$
with $x_1,\dots,x_k, x'_1,\dots,x'_k\in Q$.
Note that for $\oT_2(Q)$ is commonly called
the
additive energy of $Q$ (see, e.g. \cite{TV}).
We have
$$
    \oT_k (Q) =
        p^{2k-1} \sum_\gamma |\hat {\chi}_Q (\gamma)|^{2k} \,.
$$
The following lemma is a particular case of Lemma~4 from \cite{KS}.
\begin{lemma}\label{T_k_est}
Let $Q\subset A\subset\Z_p$, $\|\chi_A \|_{A} \le K$, $k\in\N$. Then
$$\oT_k (Q) \ge \frac{|Q|^{2k}}{|A| K^{2k-2}}.$$
\end{lemma}
In particular,
\begin{equation}\label{energy_via_Wiener}
\oT_2 (A) \ge \frac{|A|^3}{\|\chi_A\|_{A}^2}.
\end{equation}

For subsets $A,B$ of an ambient additive abelian group their sum and difference
are defined
in
a natural way:
$$A\pm B=\{a\pm b:\,a\in A, b\in B\}.$$
The following result is
the current version of the Balog-- Szemer\'edi-- Gowers
theorem \cite{Schoen_BSzG} (see also \cite{BG}).
\begin{lemma}\label{BSG} If $G$ is an additive abelian group,
$A$ is a nonempty finite subset of $G$, $\oT_2 (A)\ge|A|^3/L$, then
there exists $A'\subset A$ such that $|A'|\gg|A|/L$ and
\begin{equation}\label{differenceset}
|A'-A'|\ll L^4|A'| \,.
\end{equation}
\end{lemma}
Next, it is known that
$$|A'||A'+A'|\le|A'-A'|^2$$
(see Corollary~6.29 from \cite{TV}). Hence, (\ref{differenceset}) implies the inequality
\begin{equation}\label{sumset}
|A'+A'|\ll L^8|A'|.
\end{equation}

Another important ingredient from Additive Combinatorics is Freiman's theorem.
Define a generalized arithmetic progression (GAP) as a subset of $\Z_p$
of the form
$$P=P(x_0;\ox;\ow) = \left\{x_0+\sum_{i=1}^d v_ix_i:\,\,0\le v_i<w_i
\,(i=1,\dots,d)\right\}$$
where $\ox=(x_1,\dots,x_d)\in \Z_p^d$, $\ow=(w_1,\dots,w_d)\in\N^d$.
We will assume that all $x_i$ are not equal to zero.
The dimension of $P$ is $d$ and the size of $P$ is $\prod_{i=1}^d w_i$.
The following result is the current version of the Freiman's theorem \cite{Sanders}.
\begin{lemma}\label{Freiman} If $B$ is a nonempty subset
of $\Z_p$, $|B+B|\le M|B|$, $M\ge2$, then there is a GAP $P$ of dimension
at most $\log^{3+o(1)}M$ and size at most $|B|$ such that
$$|B\cap P|\ge|B|\exp\left(-\log^{3+o(1)}M\right).$$
\end{lemma}

Applying subsequently (\ref{energy_via_Wiener}), Lemma~\ref{differenceset}
with (\ref{sumset}), and Lemma~\ref{Freiman} we get
\begin{lemma}\label{meetGAP} For any $\ve>0$ and $K\ge K(\ve)$
if $A$ is a nonempty subset of $\Z_p$ with $\|\chi_A\|_A\le K$
and
\begin{equation}\label{d0_def}
d_\ve=d_\ve(K)=\log^{3+\ve}K
\end{equation}
then there exists a GAP $P$ of dimension at most $d_\ve$ and size at most
$|A|$ such that
$$|A\cap P|\ge|A| e^{-d_\ve}.$$
\end{lemma}

Our immediate purpose is to put some multiplicative translate of a set with
small Wiener norm into a small segment of $\Z_p$. To do it, recall
Blichtfeld's lemma (\cite{TV}, Lemma~3.27).
\begin{lemma}\label{Blichtfeld} Let $\Gamma\subset\R^d$ be a lattice
of full rank, and let $V$ be an open set in $\R^d$ such that
$\mes(V)>\mes(\R^d/\Gamma)$. Then there exist distinct $x,y\in V$
such that $x-y\in\Gamma$.
\end{lemma}

Let $P=P(x_0;\ox;\ow)$ be the GAP from Lemma~\ref{meetGAP}, let
$$\alpha_i=\frac{(|A|/p)^{1/d}}{w_i}$$
for $i=1,\dots,d$, $\delta>0$ be a small number,
$$V_\delta=\prod_{i=1}^d (-\delta,\alpha_i+\delta)\subset\R^d.$$
We observe that
$$\mes(V_\delta)>\prod_{i=1}^d \alpha_i=\frac{|A|}p\prod_{i=1}^d w_i^{-1}\ge\frac1p.$$
Let $\Gamma$ be the lattice
$$\Gamma=\Z^d+\frac{\ox}p\Z.$$
Then $\Gamma$ is a union of $p$ translates of $\Z^d$. Consequently,
$\mes(\R^d/\Gamma)=1/p$. Now we can apply Lemma~\ref{Blichtfeld}
and conclude that there exist distinct $x,y\in V_\delta$ such that
$x-y\in\Gamma$. Tending $\delta$ to $0$ we see that there are distinct
points
$$x,y\in V_0=\prod_{i=1}^d [0,\alpha_i]$$
with $x-y\in\Gamma$. Equivalently, putting
$$\Z_p^*=\Z_p\setminus\{0\}$$
and denoting
by $|z|$, $z\in \Z_p$ the minimal absolute value of
a representative of $z$ in $\Z$,
we see that
there exists $q\in\Z_p^*$,
$q<p$ such that for $i=1,\dots,d$
the following holds
$|qx_i|\le p\alpha_i$.

For any $x\in P$ we have
$$|q(x-x_0)|=\left|q\sum_{i=1}^d v_ix_i\right|
<\sum_{i=1}^d w_i|qx_i|\le\sum_{i=1}^d w_i\alpha_i=dp(|A|/p)^{1/d}.$$
So, we get the following structural property of sets with small Wiener norm.
\begin{lemma}\label{lin_struc} For any $\ve>0$ and $K\ge K(\ve)$
if $A$ is a nonempty subset of $\Z_p$ with $\|\chi_A\|_A\le K$,
$d_\ve$ is defined by (\ref{d0_def}),
$$m=\left[d_{\ve} p\left(\frac{|A|}p\right)^{1/d_\ve}\right],$$
then there exist $x_0\in\Z_p$ and $q\in\Z_p^*$ such that for the set
$$B=q(A-x_0)=\{q(x-x_0):\,x\in A\}$$
we have
$$|B\cap [-m,m]|\ge|A| e^{-d_\ve}.$$
\end{lemma}
\bigskip
\section{Upper estimates of $T_k(Q)$ for scattered $Q$}
\bigskip

Let us formulate the main result of the section.

\begin{lemma}
Let $I,k,m,M$ be positive integers.
    Let also $Q = \bigsqcup_{i=1}^I Q_i \subseteq \Z$ be a set such that
    $Q_i \subseteq [-4^{i}m, - \frac{4^i}{2}m ) \cup ( \frac{4^i}{2}m, 4^im]$,
$i$ runs over a subset of $\N$ of cardinality $I$, and $|Q_i| = M$.
    Then
    \begin{equation}\label{f:T_k_4}
        \oT_k (Q) \le 2^{8k} k^k I^k M^{2k-1} \,.
    \end{equation}
\label{l:3}
\end{lemma}

{\bf Proof of Lemma \ref{l:3}.}
First of all,  put $Q^{+} = Q\cap \{ x ~:~ x \ge 0 \}$ and $Q^{-} = Q\setminus Q^{+}$.
Using  H\"{o}lder inequality, one can easily obtain
$$
    \oT_k (Q) \le 4^k \max \{ \oT_k (Q^{+}), \oT_k (Q^{-}) \}
$$
and, thus, we need in an appropriate upper bound for $\oT_k (Q^{+}), \oT_k (Q^{-})$.
Without loosing of generality, we bound just $\oT_k (Q^{+})$, and, moreover, we write $Q$ instead of $Q^{+}$.

Further, put $N_k (x) = |\{ q_1 + \dots + q_k = x ~:~ q_j \in Q \} | \,.$
Clearly, $\sum_x N^2_k (x) = \oT_k (Q)$ and
$$
    \sum_x N_k (x) = |Q|^k = I^k M^k \,.
$$
In view of the last identity it is sufficient to prove the following uniform estimate for $N_k (x)$.

\begin{lemma}
    For any $x$, we have
$$
    N_k (x) \le 2^{6k} k^k M^{k-1} \,.
$$
\end{lemma}
{\bf Proof of the lemma.}
Take
a vector $\v{s} = (s_1, \dots,s_b)$, $s_1+\dots +s_b = k$, and put
$$
    N^{\v{s}}_k (x) = | \{ q_1+\dots + q_k = x ~:~ \exists s_1 \mbox{ elements from } A_{i_1}
        , \dots ,
                           \exists s_b \mbox{ elements from } A_{i_b} \} | \,,
$$
where  $i_1 < i_2 < \dots < i_l$.
Then
\begin{equation}\label{f:N_k_1}
    N_k (x) = \sum_{\v{s}} N^{\v{s}}_k (x) \cdot \frac{k!}{s_1! \dots s_b!} \,.
\end{equation}
Thus, we need to estimate $N^{\v{s}}_k (x)$ for any $\v{s}$.
Because of
\begin{equation}\label{f:N_k_2}
    N^{\v{s}}_k (x) \le \sum_{q_1 \in A_{i_1}} \dots \sum_{q_b \in A_{i_{b-1}}} \delta_0 (q_1+\dots+q_b - x)
        \le
            \D_1 (\v{s}) \dots \D_{b-1} (\v{s}) M^{k-1} \,,
\end{equation}
where $\D_l (\v{s})$ is the number of choices for indices of sets $A_{i_l}$,
and $\delta_0 (z)$ is the function such that $\delta_0 (z) =1$ iff $z=0$.
We need to estimate the quantities $\D_l (\v{s})$.
Suppose that the sets $A_{i_1},\dots,A_{i_{l-1}}$ are fixed and let us find an upper bound for the number of
sets $A_{i_l}$.
Let $z$ be the least integer number such that
\begin{equation}\label{f:z_def}
    \sum_{j=1}^{l-1} s_j 4^j \le s_l \frac{4^{l+z}}{2} \,.
\end{equation}
Then the number of the sets $A_{i_l}$ is bounded by $z+1$.
Indeed, without loosing of generality, we can suppose that $i_j = j$, $j \in [l-1]$
and $i_l = l + z'$, $z'>z$.
Then the set $A_{i_l}$ is defined uniquely because otherwise we have a solution of the equation
\begin{equation}\label{tmp:17.03.2014_1}
    \mu_1 + \dots + \mu_{l-1} + \mu_l = x = \mu'_1 + \dots + \mu'_{l-1} + \mu'_l \,,
\end{equation}
where $\mu_j, \mu'_j \in s_j A_{i_j}$, $j \in [l-1]$, and, similarly,  $\mu_l \in s_l A_{l+z'}$,
$\mu'_l \in s_l A_{i_l}$, $i_l < l+z'$.
If (\ref{tmp:17.03.2014_1}) takes place then
$$
    s_l \frac{4^{l+z}}{2} \le s_l \frac{4^{l+z'}}{2} < \mu'_l - \mu_l \le \mu_1 + \dots + \mu_{l-1} \le \sum_{j=1}^{l-1} s_j 4^j
$$
with a contradiction.
It follows that
$$
    \D_l (\v{s}) \le \log (2\sum_{j=1}^{l-1} s_j 4^{j-l} ) + 1
        \le
            \log ( 2\max_{1\le j \le l-1} \{ s_j 2^{j-l} \} ) + 1 \,.
$$
Let $m_1 < m_2 < \dots < m_t$ be the local maximums of the sequence
$\max_{1\le j \le l-1} \{ s_j 2^{j-l} \}$, $l \in [b-1]$.
Let also $d_j$ be the number of appearing of
the
maximum $m_j$.
Then $\sum_{j=1}^t d_j = k$.
Further, by the construction of the sequence $\max_{1\le j \le l-1} \{ s_j 2^{j-l} \}$, $l \in [b-1]$
one can see that $d_j \le \log 2s_j$, $j\in [t]$.
Returning to (\ref{f:N_k_1}), and having (\ref{f:N_k_2}), we get
$$
    N_k (x)
        \le
            M^{k-1} \sum_{\v{s}} \frac{k!}{s_1! \dots s_b!} \cdot
                (\log 2s_{m_1} + 1)^{d_1} \dots (\log 2s_{m_t} + 1)^{d_t}
                    \le
$$
$$
    \le
        M^{k-1} e^{k} k! \sum_{s_{m_1},\dots,s_{m_t}} \prod_{j=1}^t \frac{(\log 2s_{m_j} + 1)^{\log 2s_{m_j}}}{s_{m_j}!}
            \le
$$
$$
            \le
                M^{k-1} e^{2k} k! \left( \sum_s \frac{(\log 2s + 1)^{\log 2s}}{s^s} \right)^t
                    \le
                        2^{6k} k^k M^{k-1}
$$
as required.
Thus, we have proved
our
lemma and, hence, Lemma \ref{l:3}.
$\hfill\Box$

\begin{note}
    If one allows an additional multiplies of the form $(\log k)^k$ in bound (\ref{f:T_k_4})
    then the result follows immediately. Indeed, we can split our set $A$ onto sets $B_1,\dots,B_r$, $r\sim \log k$
    such that each $B_j$ contains $A_{l}$ with $l\equiv j \pmod r$.
    Thus we lose exactly $(\log k)^k$ multiple but any set $A_{i_l}$ in each $B_j$ is defined uniquely, all $\Delta_j (\v{s}) = 1$ (see formulas (\ref{f:N_k_2}), (\ref{f:z_def})),
    and, hence, $\oT_k (B_j) \le C^k k^k M^{k-1} |B_j|^k$, where $C>0$ is an absolute constant.
\end{note}

\bigskip
\section{Proof of Theorem~\ref{mediumsize}}
\bigskip
We fix an arbitrary $\ve>0$ and assume that
\begin{equation}\label{K_assum}
\|\chi_A\|_A\le K,\quad K_\ve\le K\le
(\log(p/|A|))^{1/3}(\log\log(p/|A|))^{-1-\ve}.
\end{equation}
Our aim is to prove that (\ref{K_assum}) cannot hold provided that
$p/|A|$ exceeds some quantity depending on $\ve$. Since $\ve>0$
is arbitrary, the theorem will follow.

We take $x_0, q, m,$ and $B$ accordingly with Lemma~\ref{lin_struc}.
Since
$$\hat{\chi}_B (\gamma)=e_p(-qx_0\gamma)\hat{\chi}_A (q\gamma),$$
we conclude that $\|\chi_B\|_A=\|\chi_A\|_A$. Thus,
\begin{equation}\label{B_norm}
\|\chi_B\|_A\le K.
\end{equation}

Let $l_0$ be the maximal positive integer $l$ with $2^lm<p/3$,
$$D_l=\{b\in B:\,|b|\le 2^lm\},\quad 0\le l\le l_0,$$
$$\eta=\exp(-CK)$$
for a large constant $C$, and
$$M=\left[\eta|A| e^{-d_\ve}\right].$$
If for some $l\ge1$ we have $|D_l\setminus D_{l-1}|<M$
then applying Lemma~\ref{trig} to $n=2^{l-1}m$
and taking into account the inequality $|D_l|\ge|D_0|$
and the lower bound for $|D_0|$ from Lemma~\ref{lin_struc}
we find
$$\|\hat{\chi}_B \|_1\gg\min\left(\log\frac1{\eta}, \log|D_0|\right).$$
Since
$$\log|D_0|\ge\log|A|-d_\ve\gg(\log p/\log\log p)^{1/3}
>K(\log\log p)^{2/3}>\log\frac1{\eta} \,,$$
we see that
$$\|\hat{\chi}_B \|_1\gg\log\frac1{\eta},$$
and we get contradiction with~(\ref{B_norm}) provided that $C$ is large
enough.

Thus, it is enough to consider the case where $|D_l\setminus D_{l-1}|\ge M$
for all $l=1,\dots,l_0$. For each $l$
with $l\equiv 0 \pmod 2$ we take $S_l\subset D_l\setminus D_{l-1}$
with $|S_l|=M$. Define
$$Q=\bigsqcup_l S_l.$$
Now we are in position to use Lemma~\ref{l:3} with $k=[K]$
and the sets $Q_i$ that are the sets $S_l$ in another numeration ($I=[l_0/2]$).
Let us compare the upper estimate~(\ref{f:T_k_4}) for $\oT_k (Q)$
with the lower estimate from
Lemma~\ref{T_k_est}
taking into account that
$|Q|=IM$. After simple calculations we obtain
$$\frac{|Q|}{|A|}I^{k-1}\le K^{3k-2}2^{8k}$$
implying (because of $|Q|/|A|\le\exp(\log^{3+\ve}K)$)
\begin{equation}\label{last}
I\ll K^3.
\end{equation}
We have
$$I\ge l_0/2-1\gg\log(p/m)\ge d_\ve^{-1}\log(p/|A|)-\log d_\ve \,.$$
Recalling (\ref{K_assum}) and (\ref{d0_def}) we see that
$$|I|\gg d_\ve^{-1}\log(p/|A|) \gg \log(p/|A|)(\log\log(p/|A|))^{-3-\ve}.$$
So, (\ref{last}) does not agree with (\ref{K_assum})
as required.
$\hfill\Box$
\bigskip


{}
\bigskip

\noindent{S.V.~Konyagin\\
Steklov Mathematical Institute,\\
ul. Gubkina, 8, Moscow, Russia, 119991}
\\
and
\\
MSU,\\
Leninskie Gory, Moscow, Russia, 119992\\
{\tt konyagin@mi.ras.ru}

\bigskip

\noindent{I.D.~Shkredov\\
Steklov Mathematical Institute,\\
ul. Gubkina, 8, Moscow, Russia, 119991}
\\
and
\\
IITP RAS,  \\
Bolshoy Karetny per. 19, Moscow, Russia, 127994\\
{\tt ilya.shkredov@gmail.com}


\begin{thebibliography}{99}


\bibitem{BG}
{\sc J.~Bourgain and M.Z.~Garaev. }
On a variant of sum-product estimates and explicit exponential sum bounds in prime fields //
Math. Proc. Cambridge Philos. Soc., {\bf 146}:1, 1--21, 2009.


\bibitem{GK}
{\sc B.J.~Green, S.V.~Konyagin. } On the Littlewood problem modulo a prime // Canad. J. Math., {\bf 61}:1, 141–-164, 2009.

\bibitem{Kon}
{\sc S.V.~Konyagin. } On a problem of Littlewood // Izvestiya of Russian Academy of Sciences, {\bf 45}:2, 243--265, 1981.

\bibitem{KS}
{\sc S.V.~Konyagin, I.D.~Shkredov. }
Quantitative version of Beurling--Helson theorem // submitted.

\bibitem{Lebedev}
{\sc V.V.~Lebedev. } Absolutely convergent Fourier series. An improvement of Beurling--Helson theorem. // Funct. Anal. Appl., {\bf 46}:2, 52--65, 2012.

\bibitem{MGPS}
{\sc O.C.~McGehee, L.~Pigno, B.~Smith. }
Hardy's inequality and the $L^1$ norm of exponential sums // Annals of Math., {\bf 113}, 613--618, 1981.

\bibitem{Sanders}
{\sc T.~Sanders. } The Littlewood--Gowers problem // J. Anal. Math., {\bf 101}, 123--162, 2007.

\bibitem{Sanders_add}
{\sc T.~Sanders. }
The structure theory of set addition revisited // Bull. AMS, {\bf 50}:1, 93--127, 2013.


\bibitem{Schoen_BSzG}
{\sc T.~Schoen. } New bounds in Balog--Szemer\'{e}di--Gowers theorem // Combinatorica, accepted.

\bibitem{TV}
{\sc T.~Tao, V.~Vu. }
Additive combinatorics / CUP, 2006.

\bibitem{Zyg}
{\sc A.~Zygmund. } Triginometric series / V.~2, CUP, 2002.


\end{thebibliography}
\end{document}